\documentclass[reqno, 12pt]{amsart}
\def\C{{\mathbb C}}

\def\Z{{\mathbb Z}}
\def\Zp{{\mathbb Z}_p}

\def\Fq{{{\mathbb F}_q}}

\newtheorem{lemma}{Lemma}

\newtheorem{prop}{Proposition}

\newtheorem{theorem}{Theorem}

\theoremstyle{definition}
\newtheorem{defn}{Definition}

\theoremstyle{remark}
\newtheorem{rem}{Remark}

\begin{document}
\title[On the $L$-series of F.\ Pellarin]
{On the $L$-series of F.\ Pellarin}
\author{David Goss}
\thanks{Dedicated to the memory of my dear friend David Hayes}
\address{Department of Mathematics\\The Ohio State University\\231
W.\ $18^{\rm th}$ Ave.\\Columbus, Ohio 43210}

\email{goss@math.ohio-state.edu}

\date{July 24, 2011}

\begin{abstract}
The calculation, by L.\ Euler, of the values at positive even integers of the Riemann
zeta function, in terms of powers of $\pi$ and rational numbers, 
was a watershed event in the history of number theory and classical analysis. 
Since then many important analogs involving $L$-values and periods
have been obtained. In analysis in
finite characteristic, a version of  Euler's result
was given by L.\ Carlitz \cite{ca2} in the 1930's which involved
the period of a rank $1$ Drinfeld module (the Carlitz module) in place
of $\pi$. In a very original work \cite{pe2},
F.\ Pellarin has quite recently established a ``deformation'' of Carlitz's
result involving certain $L$-series and the deformation of the Carlitz period
given in \cite{at1}. Pellarin works only with the values of this $L$-series 
at positive integral points. We show here how the techniques of \cite{go1} also 
allow these new $L$-series to be analytically continued -- with associated
trivial zeroes -- and interpolated at finite primes.
\end{abstract}

\maketitle

\section{Introduction}\label{intro} 
In 1734, after a number of attempts, Euler succeeded in giving a closed form
formula for $\zeta(2)$ where $\zeta(s)=\sum_{n=0}^\infty n^{-s}$ is the Riemann zeta 
function. Indeed, at that time Euler
obtained the famous formula $\zeta(2)=\pi^2/6$ (as well as
$\zeta(4)=\pi^4/90$), see e.g., \cite{ay1}. This marks the beginning of the remarkably 
profound association between periods and $L$-series of classical motives; a subject 
being vigorously investigated to this day.

Let $p$ be a prime and let $q=p^{m_0}$ for a fixed positive integer $m_0$. Set
$A:=\Fq[\theta]$ where $\theta$ is an indeterminate. Approximately 200 years
after Euler, in the mid-1930's, L.\ Carlitz \cite{ca2} used characteristic
$p$ analysis to establish an analog for $A$ of
Euler's result.  More precisely, 
let $k:=\Fq(\theta)$ and $K:=\Fq((1/\theta))$ (which is the completion of 
$k$ at the infinite place). Let $j$ be a positive integer and put
\begin{equation}\label{intro1}
\zeta_A(j):=\sum_{f(\theta)~\rm monic} f^{-j}\,;\end{equation}
 this series is readily seen to converge to  nonzero element of $K$. Carlitz then 
established the existence of a nonzero constant $\tilde{\pi}$ so that if
$(q-1)\mid j$, then $\zeta(j)/\tilde{\pi}^j\in k$. In the 1970's, this result
was independently rediscovered by the present author in the context of
Eisenstein series.

The constant $\tilde{\pi}$, which is uniquely defined up to an element of
$\Fq^\ast$, is the period of the {\it Carlitz module} $C$ \cite{ca1}; i.e., 
$\tilde{\pi}$ is generator of the rank $1$ $A$-module of periods 
associated to the exponential function $e_C(z)$ of $C$. (As such, $\tilde{\pi}$
is actually analogous to $2\pi i$ in classical theory.)

Carlitz also gave a beautiful product formula for his period which was remarkably
simplified in  \cite{at1} (see also \cite{abp1})  in the following fashion. Let $\theta_1$ be a 
fixed choice of $(q-1)$-st root of $-\theta$. Then one has
\begin{equation}\label{intro2}
\tilde{\pi}=\theta_1 \theta \prod_{i=1}^\infty \left( 1-\theta^{1-q^i}\right)^{-1}\,.
\end{equation}
In \cite{at1} (and \cite{abp1}),
this formula is also ``deformed'' in the following fundamental
fashion. Let $t$ be another indeterminate and set
\begin{equation}\label{intro3}
\Omega(t):=\theta_1^{-q}\prod_{i=1}^\infty \left(1-t/\theta^{q^i}\right)\,.\end{equation}
It follows immediately that the product for $\Omega(t)$ converges
for all $t\in K$, and one has 
\begin{equation}\label{intro3.1}
\Omega(\theta)\tilde{\pi} =-1\,.
\end{equation}

The function $\Omega(t)$ is essential for the study of tensor powers of the
Carlitz module and special $\Gamma$-values. As such, it is also
natural to search for a 
connection between it and characteristic $p$ $L$-series. This has been 
very recently accomplished
in seminal work by Federico Pellarin. Indeed, in \cite{pe2} (which
depends of \cite{pe1}),  Pellarin obtains the
following elegant generalization of the formulas of Euler and Carlitz. Let 
$t\in K$ and let $\chi_t$ be the quasi-character of $A$ given simply by
$$\chi_t(f):=f(t)\,,$$
for $f(\theta)\in A$. Let $j$ be a positive integer and set
\begin{equation}\label{intro4}
L(\chi_t,j):=\sum_{f(\theta)~\rm monic} \chi_t(f)/f^j\,.\end{equation}
It is easy to see that this converges for $t$ sufficiently small. 

Let $j$ now be congruent to $1$ modulo $q-1$. 
Then, also
using Eisenstein series (or, rather, vectorial Eisenstein series), 
Pellarin shows the existence of a rational function
$b_j(\theta,t)\in \Fq(\theta,t)$ such that
\begin{equation}\label{intro5}
L(\chi_t,j)=b_j(\theta,t)\tilde{\pi}^j\Omega(t)\,.\end{equation}
(N.B.: The normalized function $\tilde{\pi}\Omega(t)$ is independent of the
choice of $\theta_1$.)  As $\tilde{\pi}\Omega(\theta)=-1$, Equation \ref{intro5} is
clearly a deformation of the result of Carlitz.

One can compute that $b_1(\theta,t)\equiv -1$. As such,
\begin{equation}\label{intro6}
\lim_{t\to \theta} L(\chi_t,1)=1
\end{equation}
in agreement with previous results on zeta values.

In his paper, Pellarin only considers positive values of $j$. Clearly, however, one
would like these ``Pellarin $L$-series'' to have the good analytic properties possessed by
the $L$-series of Drinfeld modules. It is our purpose here to establish these
properties. Indeed, we show that the needed results follow from the elementary
estimates of Lemma 8.8.1 in \cite{go1}.

The $L$-series considered here are but the first in a long line of functions one now
wants to understand. In fact, the quasi-character $\chi_t$ may be used to deform
any $L$-series of any Drinfeld module, $A$-module, etc., (see e.g., 
\cite{bo1}, \cite{bo2}, \cite{go2}) and all such functions must
be studied. 

Moreover, $A=\Fq[\theta]$ is but the simplest base ring in the general theory
of Drinfeld modules. One also expects similar results for general $A$ where there
are obviously some subtleties involved. Still any such formula will involve
the periods of the Hayes modules \cite{ha1}; these objects
are normalized rank $1$ Drinfeld
$A$-modules which were constructed by David Hayes as the correct
generalizations of the Carlitz module. Indeed for $A=\Fq[\theta]$,
the Hayes module is the Carlitz module.
Therefore, it is only fitting that Hayes' ideas appear in the present context
as his work has been absolutely
essential in moving
the arithmetic theory from the simplest case of $\Fq[\theta]$ to Drinfeld's arbitrary base
rings $A$. In fact, without his construction in the rank $1$ case almost no
arithmetic in general is even possible.
As such, it is my pleasure to dedicate this paper to David's memory
and mathematical legacy.

I also express my gratitude to Federico Pellarin for his remarks on a previous version
of this work.

\section{Definition of the Pellarin $L$-series}\label{definitions}
As before, let $A=\Fq[\theta]$, $q=p^{m_0}$; we put $k:=\Fq(\theta)$,
$K:=k_\infty=\Fq((1/\theta))$. Let $A_+$ be the set of monic elements of $A$.
For $d$ a nonnegative integer, we let $A_+(d)$ be the monics of degree $d$ in $A$
and $A(d)$ the vector space of all polynomials of degree $\leq d$.
As usual, $K$ comes equipped with the absolute value $\vert?\vert_\infty$ such that $\vert \theta\vert_\infty=q$.  We
let $\overline{K}$ be a fixed algebraic closure of $K$ equipped with the canonical
extension of $\vert?\vert_\infty$ and we let $\C_\infty$ be its completion. Let
$v_\infty(?)$ be the associated canonical additive valuation on $\C_\infty$ extending
the valuation $v_\infty$ on $K$ with $v_\infty(1/\theta)=1$.

Let $\pi$ be a fixed uniformizer in $K$ of the form $\pi=1/\theta+
\{\rm higher~terms\}$. (In fact, as will be seen in Remark \ref{plseries0.000001} just
below, the reader may take $\pi=1/\theta$ without
any loss of generality.) Let $a\in A_+$ and put
\begin{equation}\label{plseries.00001}
\langle a\rangle_\pi:=\pi^{\deg a}a\,.
\end{equation}
\noindent
Clearly, the map $a\mapsto \langle a \rangle_\pi$ is a homomorphism from the
multiplicative monoid $A_+$ into the
group of $1$-units of $K$. The Binomial Theorem therefore allows one to
raise $\langle a\rangle_\pi$ to any power $y\in \Z_p$.

\begin{defn}\label{plseries.001}
We set
\begin{equation}\label{plseries0}
{\mathbb S}_\infty:=\C_\infty^\ast \times \Z_p\,.
\end{equation}
\end{defn}
The space $\mathbb S_\infty$ is a commutative group whose operation will be written
additively. Let $a\in A_+$ and $s=(x,y)\in \mathbb S_\infty$. Then we define
\begin{equation}\label{plseries0.1}
a^s:=x^{\deg a}\langle a \rangle_\pi^y \,.
\end{equation}
\noindent
Let $j$ be an integer and set $s_j:=(\pi^{-j},j)\in {\mathbb S}_\infty$. One then
readily checks that $a^{s_j}=a^j$ where ``$a^j$'' has the usual meaning.
We will freely refer to $s_j$ as ``$j$''. 

The space $\mathbb S_\infty$ is the domain of characteristic $p$ $L$-series at the place
$\infty$ of $k$.

\begin{rem}\label{plseries0.000001}
Let $\pi_1$ and $\pi_2$ be as above. Note that for $a\in A_+$ of degree $d$, and
$y\in \Zp$,  we
have $\langle a \rangle_{\pi_1}^y =(\pi_1/\pi_2)^{dy}\langle a \rangle_{\pi_2}^y$.
Clearly $\vert (\pi_1/\pi_2)^{dy}\vert_\infty=1$. Therefore, as we are interested here
in establishing certain estimates, this factor is harmless. {\em As such, we now
set $\pi:=1/\theta$ for the rest of this paper and drop the reference to $\pi$.}\end{rem}

\begin{defn}\label{plseries1}
Let $t\in \C_\infty$ and $f\in A$. Then we define $\chi_t(f):=f(t)$. 
\end{defn}
\noindent
Obviously, $\chi_t$ is just the evaluation map of $f(\theta)$ at $t\in \C_\infty$, 
$\theta \mapsto t$,  and thus
is clearly an $\Fq$-algebra morphism from $A$ to $\C_\infty$. We view it as
a quasi-character of $A$.

\begin{defn}\label{plseries2}
For $\beta$ a nonnegative integer and $s\in \mathbb S_\infty$, we formally define
\begin{equation}\label{plseries3}
L(\chi_t^\beta,s):=\sum_{a\in A_+}\chi_t(a)^\beta a^{-s}=\prod_{f~\rm monic~prime}
(1-\chi_t^\beta(f) f^{-s})^{-1}\,.\end{equation}
We call $L(\chi_t^\beta,s)$ the {\it Pellarin $L$-series associated to $\chi_t^\beta$}.
\end{defn}

\begin{lemma}\label{plseries4}
Let $t\in \C_\infty$ and set $\lambda_t:=\max\{1,\vert t\vert_\infty\}$. Let $f(\theta)\in
A$ have degree $d$. Then $\vert f(t)\vert_\infty\leq \lambda_t^d$.
\end{lemma}
\begin{proof} Obvious.\end{proof}

\begin{prop}\label{plseries5} Let $t\in \C_\infty$ and $\lambda_t$ be as
{\rm Lemma \ref{plseries4}}.
Then the Euler-product for the $L$-series $L(\chi_t^\beta,s)$ converges on the
``half-plane'' of
$\mathbb S_\infty$ given by $\{s=(x,y)\mid \vert x\vert_\infty>\lambda_t^\beta\}$.
\end{prop}
\begin{proof} This follows immediately  from the definitions and
 Lemma \ref{plseries4}.\end{proof}

In Theorem \ref{main6} we establish that $L(\chi_t^\beta,s)$ can be analytically continued to an
entire function on $\mathbb S_\infty$ in the sense of Section 8.5 of \cite{go1}.
\section{Review of estimates}\label{estimates} 
We recall here the basic estimates from  Section 8.8 of \cite{go1} necessary for us
and refer the reader there for
their derivations. Let $J_0$ and $J_1$ be two fields over $\Fq$.  Let $W\subseteq J_0$ be a 
finite dimensional $\Fq$-vector space of dimension $d$. Further let
$\{{\mathcal L}_1,\dots ,{\mathcal L}_t\}$ be $\Fq$-linear maps of $J_0$ into $J_1$. Let
$x\in J_0$ and let $\{i_1,\dots,i_t\}$ be nonnegative integers so that
\begin{equation}\label{est1}
\sum_{h=1}^t i_h<(q-1)d\,.
\end{equation} 

\begin{lemma}\label{est2}
Under the above assumptions we have 
\begin{equation}\label{est3}
\sum_{w\in W}\left(\prod^t_{h=1}{\mathcal L}_h(x+w)^{i_h}\right)
=0\,.
\end{equation}\end{lemma}

Assume now that $J_1$ has an additive discrete valuation $v$ with
$v({\mathcal L}_h(w))>0$ for all $h$ and $w$. Let $\{ i_h\}$ now be an arbitrary collection of 
non-negative integers, and for $j>0$ put
$$W_j:=\{ w\in W\mid v({\mathcal L}_h(w))\ge j~~{\rm for~all~}h\}\,.$$
Finally, set $Q:=\sum_j \dim_{\Fq}W_j$. 
\begin{lemma}\label{est4}
Under the above assumptions, we have
\begin{equation}\label{est5}
v\left(\sum_{w\in W}\prod^t_{h=1}{\mathcal L}_h(w)^{i_h}\right)\ge 
(q-1)Q\,.\end{equation}\end{lemma}

\section{The main theorem}\label{main}
We now establish the analytic continuation of $L(\chi_t^\beta,s)$ using
Lemma \ref{est4}. We note that it is not always possible to apply this result
directly and so we proceed in a slightly indirect fashion. 

The first step is to rewrite the Pellarin $L$-series in the usual fashion by summing
according to degrees. That is, upon unraveling the definitions,
we find in the half-plane of convergence that
\begin{equation}\label{main0}
L(\chi_t^\beta,s)=\sum_{j=0}^\infty x^{-j}\left(\sum_{a\in A_+(j)} \chi_t(a)^\beta 
\langle a \rangle^{-y}\right)\,.\end{equation} 
For each fixed $y\in \Z_p$, the goal is to establish that $L(\chi_t,x,y)$ is an entire
power series in $x^{-1}$ with the resulting function on $\mathbb S_\infty$ having
good continuity properties. 

Let $t\in \C_\infty$ and set $\delta_t:=\max\{-v_\infty(t)+1,1\}$. Let
$0\neq \alpha\in \C_\infty$ be chosen so that $-v_\infty(\alpha)\geq\delta_t$.
\begin{lemma}\label{main1}  Let $\{t,\alpha\}$ be as above.
Let $d$ be a nonnegative integer and $0\leq i\leq d$.
Then, under the above assumption on $\alpha$, we have
\begin{equation}\label{main2}
v_\infty(t^i/\alpha^d)\geq d-i\,.\end{equation}\end{lemma}
\begin{proof}
Obviously one has
\begin{equation}\label{main3}
v_\infty(t^i/\alpha^d)=iv_\infty(t)-v_\infty(\alpha)d\,.
\end{equation}
There are then two cases: $v_\infty(t)\geq 0$ and $v_\infty(t)<0$. In the
first case, by assumption, $-v_\infty(\alpha)\geq 1$ and the result follows directly.
In the second case, we have $-v_\infty(\alpha)\geq -v_\infty(t)+1$. Thus,
$$v_\infty(t^i/\alpha^d)=iv_\infty(t)-v_\infty(\alpha)d\geq (-v_\infty(t)+1)d+iv_\infty
(t)=d+v_\infty(t)(i-d)\geq d-i\,.$$
\end{proof}
With the above choices of $\{t,\alpha\}$  we now temporarily 
set 
\begin{equation}\label{main4}
L_\alpha(\chi_t^\beta,s):=L(\chi_t^\beta,\alpha^\beta x,y)=
\sum_{j=0}^\infty x^{-j}\left(\sum_{a\in A_+(j)} (\chi_t(a)/\alpha^j)^\beta 
\langle a \rangle^{-y}\right)\,.\end{equation}

Let $\Fq[1/\theta](j)$ the $\Fq$-vector space of polynomials in $1/\theta$ of degree
at most $j$. Note that if $a\in A_+(j)$ then $\langle a \rangle\in \Fq[1/\theta](j)$.
We now define two $\Fq$-linear maps $\{{\mathcal L}_1,
{\mathcal L}_2\}$ from $\Fq[1/\theta](j)$ to $\C_\infty$ as follows:
${\mathcal L}_1$ is simply the identity map, while
\begin{equation}\label{main5}
{\mathcal L}_2(\sum_{n=0}^j c_n \theta^{-n}):=\sum_{n=0}^jc_nt^{j-n}/\alpha^j\,.
\end{equation}

We can now establish our main result.
\begin{theorem}\label{main6}
The Pellarin $L$-series $L(\chi_t^\beta,s)$ analytically continues to an entire
function on ${\mathbb S}_\infty$.
\end{theorem}
\begin{proof} It is clearly sufficient to establish that $L_\alpha (\chi_t^\beta,s)$
is entire.

Let $W(j):=\{f\in \Fq[1/\theta](j)\mid f(0)=0\}$ so that $\dim_\Fq W(j)=j$.
Let $w\in W(j)$ with $v_\infty(w)=i$.
The auxiliary constant $\alpha$ has been chosen {\em precisely} to guarantee
${\mathcal L}_2(w)\geq i$ via Lemma \ref{main1}. 

Let $y\in \Z_p$ be fixed and let $c_j(y)$ be the coefficient of $x^{-j}$ in
the expansion of $L_\alpha(\chi_t^\beta,x,y)$. Lemma \ref{est4} now immediately
implies that $v_\infty(c_j(y))\geq (q-1)j(j+1)/2$. This quadratic growth is 
easily seen to be sufficient to establish the result. \end{proof}

\begin{rem}\label{main7} One also readily deduces that $L(\chi_t,s)$ is continuous
in $t$. \end{rem}

\section{Special polynomials and trivial zeroes} \label{special}
In this section, we examine the behavior of $L(\chi_t^\beta,s)$ when
$s=(x,y)$ and $y\in \Z_p$ is a nonpositive integer. So let $y=-j$, for $j\geq 0$ .
\begin{defn}\label{special1}
We set
\begin{equation}\label{special2}
z(\chi_t^\beta,x,-j):=L(\chi_t^\beta,x\pi^j,-j)=\sum_{e=0}^\infty
x^{-e}\left(\sum_{a\in A_+(e)} \chi_t^\beta(a) a^j\right)\,.
\end{equation}\end{defn}
Our next result will show that for each such $j$, $z(\chi_t^\beta,x,-j)$ is
a {\em polynomial} in $x^{-1}$ called the {\it special polynomial} at $-j$.

\begin{theorem}\label{special3}
Let $j$ be as above. Then $z(\chi_t^\beta,x,-j)\in A[t][x^{-1}]$.
\end{theorem}
\begin{proof}
As in the proof of Theorem \ref{main6}, we need to define two $\Fq$-linear
maps ${\mathcal L}_1$ and ${\mathcal L}_2$ on $A(e)$. Also as before, we set
${\mathcal L}_1$ to be the identity, and now we set ${\mathcal L}_2:=\chi_t$. Notice that
$A_+(e)=\{\theta^e+h\}$ where $h$ runs over $A(e-1)$ and
$\dim_\Fq A(e-1)=e$. Thus, Lemma \ref{est2}
immediately implies that the coefficient of $x^{-e}$ in $z(\chi_t^\beta,x,-j)$ 
vanishes when $e>(\beta+j)/(q-1)$. \end{proof}

\subsection{Trivial zeroes}\label{trivial}
Let $t=\zeta\in \C_\infty$ be a root of unity. In this case, as was pointed out
in \cite{pe1}, the quasi-character
$\chi_\zeta$ is actually a character factoring through $A/(p(\theta))$ where
$p(\theta)$ is the minimal polynomial associated to $\zeta$. As such,
the techniques of Section 8.13 of \cite{go1} are easily altered to equip
$L(\chi_\zeta^\beta,s)$ with a trivial zero at $s=-\lambda$ where $\lambda>\beta$ is a positive
integer such that $\lambda\equiv -\beta\pmod{q-1}$.

\begin{prop}\label{trivial1}
Let $\lambda>\beta$ be a positive integer congruent to $-\beta \pmod{q-1}$. Then
$L(\chi_t^\beta,s)=0$ for $s=-\lambda=(\pi^\lambda,-\lambda)\in {\mathbb S}_\infty$.
\end{prop}
\begin{proof} Theorem \ref{special3} immediately implies that 
$z(\chi_t^\beta,1,-\lambda)$ is an element of $A[t]$. On the other hand,
this element vanishes for all $t=\zeta$, where $\zeta$ is a root of unity as above.
Thus it is identically 0.\end{proof}
\begin{rem}\label{trivial2}
A combinatorial proof of Proposition \ref{trivial1} has also been given by R.\ Perkins.
\end{rem} 

\section{The $\mathfrak v$-adic theory}\label{vadic}
Let $\mathfrak v$ be a nontrivial prime of $A$ of degree $d$
and let $\C_\mathfrak v$ be the $\mathfrak v$-adic
version of $\C_\infty$. We very briefly indicate here how the above techniques may
be altered to give the $\mathfrak v$-adic interpolation of the above special polynomials
to entire functions on the $\mathfrak v$-adic analog of ${\mathbb S}_\infty$. 

Let $P(\theta)$ be the monic generator of $\mathfrak v$ and let $t\in \C_\mathfrak v$, so
that we have the $\mathfrak v$-adic quasi-character $\chi_t$ of $A$.  We now
choose the auxiliary constant $\alpha \in \C_\mathfrak v$ small enough so that
$v_\mathfrak v (\alpha^{-d} P(t))> 1$. With this choice the results now follow
directly (see also Section 8.9 of \cite{go1}).

\end{document}